\documentclass[a4paper, oneside, 11pt]{article}%
\usepackage{amsmath}
\usepackage{amsfonts}
\usepackage{amssymb}
\usepackage{graphicx}
\usepackage{appendix}
\usepackage[square, numbers, sort&compress]{natbib}
\usepackage{color}%
\providecommand{\U}[1]{\protect\rule{.1in}{.1in}}

\newtheorem{theorem}{Theorem}[section]

\newtheorem{assumption}{Assumption}

\newtheorem{lemma}[theorem]{Lemma}

\newtheorem{proposition}[theorem]{Proposition}
\newtheorem{remark}[theorem]{Remark}

\numberwithin{equation}{section}

\newcommand{\E}{{\mathbb E}}
\newcommand{\R}{{\mathbb R}}
\newcommand{\pf}{\noindent\textbf{Proof:} }
\newcommand{\eof}{\hfill{$\Box$}}

\newcommand{\BMO}{L^{2}_{\mathcal F^{W}}(0, N;\mathbb{R}^{n})}
\usepackage{url}
\usepackage{hyperref}
\hypersetup{
colorlinks=true, 
breaklinks=true, 
urlcolor=green , 
linkcolor=red, 
citecolor=blue, 
pdftitle={}, 
pdfauthor={}, 
pdfsubject={} 
}

\newcommand{\cM}{\ensuremath{\mathcal{M}}}

\usepackage[normalem]{ulem}

\begin{document}

\title{Constrained stochastic LQ control on infinite time horizon with regime switching}
\author{Ying Hu \thanks{Univ Rennes, CNRS, IRMAR-UMR 6625, F-35000 Rennes, France. Partially supported by Lebesgue
Center of Mathematics \textquotedblleft Investissements d'avenir\textquotedblright program-ANR-11-LABX-0020-01, ANR CAESARS
(No.~15-CE05-0024) and ANR MFG (No.~16-CE40-0015-01). Email:
\texttt{ying.hu@univ-rennes1.fr }}
\and Xiaomin Shi\thanks{Corresponding author. School of Mathematics and
Quantitative Economics, Shandong University of Finance and Economics, Jinan
250100, China. Partially supported by NSFC (No.~11801315), NSF of Shandong Province (No.~ZR2018QA001, ~ZR2020MA032), and the Colleges
and Universities Youth Innovation Technology Program of Shandong Province
(No. 2019KJI011). Email: \texttt{shixm@mail.sdu.edu.cn}}
\and Zuo Quan Xu\thanks{Department of Applied Mathematics, The Hong Kong Polytechnic University, Kowloon, Hong Kong.
Partially supported by NSFC (No.~11971409), Hong Kong
GRF (No.~15204216 and No.~15202817), The PolyU-SDU Joint Research Center on Financial Mathematics and the CAS AMSS-PolyU Joint Laboratory of Applied Mathematics, The Hong Kong Polytechnic University. Email: \texttt{maxu@polyu.edu.hk}}}
\maketitle

This paper is concerned with a stochastic linear-quadratic (LQ) optimal control problem on infinite time horizon, with regime switching, random coefficients, and cone control constraint. To tackle the problem, two new extended stochastic Riccati equations (ESREs) on infinite time horizon are introduced. The existence of the nonnegative solutions, in both standard and singular cases, is proved through a sequence of ESREs on finite time horizon. Based on this result and some approximation techniques, we obtain the optimal state feedback control and optimal value for the stochastic LQ problem explicitly. Finally, we apply these results to solve a lifetime portfolio selection problem of tracking a given wealth level with regime switching and portfolio constraint.



\addcontentsline{toc}{section}{\hspace*{1.8em}Abstract}

\section{Introduction}
Linear-quadratic (LQ) optimal control is one of the most important research topics in control theory. On one hand, it usually admits elegant optimal state feedback control and optimal value through the famous Riccati equation. On the other hand, it has widely applications in many fields, such as engineering, management science and mathematical finance.

Since the pioneering work of Wonham \cite{Wo}, stochastic LQ problem has been extensively studied by numerous researchers. For instance, Bismut \cite{Bi} was the first one that studied stochastic LQ problems with random coefficients. But he could not solve the related stochastic Riccati equation (SRE) in general. Tang \cite{Ta} proved the existence and uniqueness of the SRE and solved the general stochastic LQ problem with uniformly positive control weighting matrix. Chen, Li and Zhou \cite{CLZ} studied the indefinite stochastic LQ problem which is different obviously from its deterministic counterpart. Kohlmann and Zhou \cite{KZ} established the relationship between stochastic LQ problems and backward stochastic differential equations (BSDEs). Hu and Zhou \cite{HZ} solved the stochastic LQ problem with cone control constraint using Tanaka's formula. Please refer to Chapter 6 in Yong and Zhou \cite{YZ} for a systematic accounts on this subject.

Stochastic LQ problems on infinite time horizon were studied in Ait Rami and Zhou \cite{RZ}, Yao, Zhang and Zhou \cite{YZZ} via algebraic Riccati equations, linear matrix inequality and semidefinite programming techniques. While its application in tracking a financial benchmark was studied in Yao, Zhang and Zhou \cite{YZZ06}.
Sun and Yong \cite{SY} proved the equivalence of open-loop solvabilities, closed-loop solvabilities of the infinite horizon LQ problem and the existence of a static stabilizing solution to the associated algebraic Riccati equations.
Chen and Zhou \cite{CZ} addressed the conic stabilizability of the controlled stochastic differential equations with cone constraints, and solved the corresponding stochastic LQ problem on infinite horizon via stabilizing solutions of two related ESREs. Li, Zhou and Ait Rami \cite{LZR} studied a stochastic LQ problem with Markovian jumps on infinite time horizon. And for elaborate illustrations on regime switching models, one is referred to Yin and Zhang \cite{YZQ}.

All the above results on stochastic LQ problems on infinite horizon were established under the assumption that the coefficients in the problem are constants (matrices). When the coefficients are random, Guatteri and Tessitore \cite{GT} studied stochastic LQ problem on infinite horizon and infinite dimensional state space, but the control variable is absent in the diffusion of the state. Guatteri and Masiero \cite{GM} studied infinite horizon and ergodic stochastic LQ problems. Pu and Zhang \cite{PZ} studied stochastic LQ problem on infinite horizon with cone control constraint. The cost weighting matrices of the control are assumed to be uniformly positive and there is no regime switching in the above three papers.

In this paper, we study a stochastic LQ problem on infinite time horizon with regime switching and random coefficients, where the control variable has to be constrained in a cone. And the control weighting matrix in the cost functional is allowed to be possibly singular. The state process is shown to satisfy both the {$L^2$-stable} condition (see e.g. Definition 2.1 in \cite{YZZ06}) and the mean-square stabilizing condition (see e.g. Definition 2.1 in \cite{CZ} or \cite{SY}). To solve the control problem, we introduce two systems of BSDEs termed extended stochastic Riccati equations (ESREs) on infinite time horizon. Approximated by a sequence of BSDEs on finite time horizon, we prove that the two systems of ESREs admit nonnegative solutions. Eventually we solve the stochastic LQ problem with the help of the ESREs and some approximation techniques. The uniqueness of solutions for the ESREs is deduced by a verification argument.

The main contributions of this paper include at least the following three parts. Firstly, we prove that the two systems of ESREs on infinite horizon have nonnegative solutions from a point of view of BSDE without using the value function of stochastic LQ problem. Thus our method is interesting in its own right in BSDE theory and may have applications in solving other BSDEs. Secondly, in obtaining the optimal state feedback control and optimal value of the stochastic LQ problem, we do not use the stochastic maximum principle as in \cite{GM, PZ}. So our method is more directly than \cite{GM,PZ}.
Thirdly, we can handle the systems of ESREs on infinite horizon and stochastic LQ problem when the control weighting matrix is singular. To the best of our knowledge, this is the first paper concerning stochastic LQ problem on infinite time horizon with random coefficients and singular control weighting matrices.

This paper is organised as follows. In Section 2, we formulate a stochastic LQ problem on infinite time horizon with regime switching, random coefficients, and portfolio constraint. Section 3 introduces two systems of ESREs and gives some remarks about the two equations. Section 4 and Section 5 are concerned about the global solvability of two systems of ESREs and stochastic LQ problems in standard case and singular case, respectively. In Section 6, we apply the general results to solve a lifetime portfolio selection problem of tracking a given wealth level with regime switching and portfolio constraint. Some concluding remarks are given in Section 7.

\section{Problem formulation}
Let $(\Omega, \mathcal F, \mathbb{P})$ be a fixed complete probability space on which are defined a standard $n$-dimensional Brownian motion $W(t)=(W_1(t), \ldots, W_n(t))'$ and a continuous-time stationary Markov chain $\alpha_t$ valued in a finite state space $\mathcal M=\{1, 2, \ldots, \ell\}$ with $\ell>1$. We assume $W(t)$ and $\alpha_t$ are independent processes. The Markov chain has a generator $Q=(q_{ij})_{\ell\times \ell}$ with $q_{ij}\geq 0$ for $i\neq j$ and $\sum_{j=1}^{\ell}q_{ij}=0$ for every $i\in\mathcal{M}$.
Define the filtrations $\mathcal F_t=\sigma\{W(s), \alpha_s: 0\leq s\leq t\}\bigvee\mathcal{N}$ and $\mathcal F^W_t=\sigma\{W(s): 0\leq s\leq t\}\bigvee\mathcal{N}$, where $\mathcal{N}$ is the totality of all the $\mathbb{P}$-null sets of $\mathcal{F}$.

\subsection*{Notation}
We use the following notation throughout the paper:
\begin{align*}
L^2_{\mathcal{F}_t}(\Omega;\mathbb{R})&=\Big\{\xi:\Omega\rightarrow
\mathbb{R}\;\Big|\;\xi\mbox { is }\mathcal{F}_{t}\mbox{-measurable, and }\E\big(|\xi|^{2}\big)%
<\infty\Big\},\\
L^{2}_{\mathcal F}(0, \infty;\mathbb{R})&=\Big\{\phi:[0, \infty)\times\Omega\rightarrow
\mathbb{R}\;\Big|\;\phi(\cdot)\mbox{ is an }\{\mathcal{F}%
_{t}\}_{t\geq0}\mbox{-predictable process }\\
&\qquad\mbox{ \ \ \ \ with }\E\int_{0}^{\infty}|\phi(t)|^{2}dt<\infty
\Big\}, \\
L^{2, \:\mathrm{loc}}_{\mathcal F}(0, \infty;\mathbb{R})&=\Big\{\phi:[0, \infty)\times\Omega\rightarrow
\mathbb{R}\;\Big|\;\phi(\cdot)\mbox{ is an }\{\mathcal{F}%
_{t}\}_{t\geq0}\mbox{-predictable process}\\
&\qquad\mbox{ \ \ \ \ with }\E\int_{0}^{T}|\phi(t)|^{2}dt<\infty \ \mbox{for any} \ T>0
\Big\}, \\
L^{\infty}_{\mathcal{F}}(0, \infty;\mathbb{R})&=\Big\{\phi:[0, \infty)\times\Omega
\rightarrow\mathbb{R}\;\Big|\;\phi(\cdot)\mbox{ is an }\{\mathcal{F}%
_{t}\}_{t\geq0}\mbox{-predictable essentially}\\
&\qquad\mbox{ \ \ \ \ bounded process} \Big\}.
\end{align*}
These definitions are generalized in the obvious way to the cases that $\mathcal{F}$ is replaced by $\mathcal{F}^W$ and $\mathbb{R}$ by $\mathbb{R}^n$, $\mathbb{R}^{n\times m}$ or $\mathbb{S}^m$, where $\mathbb{S}^m$ is the set of symmetric $m\times m$ real matrices.

In our argument, ``almost surely" and ``almost everywhere" (a.e.) may be suppressed for notation simplicity in some circumstances when no confusion occurs.

We now introduce the problem that will be investigated in this paper.
Consider the following controlled $\mathbb{R}$-valued linear stochastic differential equation (SDE):
\begin{align}
\label{state}
\begin{cases}
dX(t)=\left[A(t, \alpha_t)X(t)+B(t, \alpha_t)'u(t)\right]dt
+\left[C(t, \alpha_t)X(t)+D(t, \alpha_t)u(t)\right]'dW(t), \ t\geq0, \\
X(0)=x, \ \alpha_0=i_0,
\end{cases}
\end{align}
where $A(t, \omega, i), \ B(t, \omega, i), \ C(t, \omega, i), \ D(t, \omega, i)$ are all $\{\mathcal{F}^W_t\}_{t\geq 0}$-predictable processes of suitable sizes for $i\in\cM$, $x\in\mathbb{R}$ and $i_{0}\in\cM$ are known, and $u(\cdot)$ is the control.
Let $\Gamma\subset\mathbb{R}^m$ be a given closed cone, i.e., $\Gamma$ is closed, and if $u\in\Gamma$, then $\lambda u\in\Gamma$, for all $\lambda\geq 0$. It represents the constraint set for controls.
The class of admissible controls is defined as the set
\begin{align*}
\mathcal{U}:=\Big\{u(\cdot)\in L^2_\mathcal{F}(0, \infty;\mathbb{R}^m)\;\Big|\; u(\cdot) \in\Gamma\Big\}.
\end{align*}
If $u(\cdot)\in\mathcal{U}$ and $X(\cdot)$ is the associated (unique) solution of \eqref{state}, then we refer to $(X(\cdot), u(\cdot))$ as an admissible pair.

Let us now state our stochastic linear quadratic optimal control problem (stochastic LQ problem, for short) as follows:
\begin{align}
\begin{cases}
\mathrm{Minimize} &\ J_\infty(x, i_0, u(\cdot))\\
\mbox{subject to} &\ (X(\cdot), u(\cdot)) \mbox{ admissible for} \ \eqref{state},
\end{cases}
\label{LQ}%
\end{align}
where the cost functional is given as the following quadratic form
\begin{align}\label{costfunctional}
J_\infty(x, i_0, u(\cdot)):=\mathbb{E}\int_0^\infty\Big(Q(t, \alpha_t)X(t)^2+u(t)'R(t, \alpha_t)u(t)\Big)dt.
\end{align}
The associated value function is defined as
\begin{align*}
V(x,i_0):=\inf_{u\in\mathcal{U}}J_\infty(x, i_0, u(\cdot)), \ x\in\mathbb{R}, \ i_0\in\cM.
\end{align*}

Throughout this paper, we put the following assumptions on the coefficients.
\begin{assumption} \label{assump1}
For all $i\in\cM$,
\begin{align*}
\begin{cases}
A(t, \omega, i)\in L_{\mathcal{F}^W}^\infty(0, \infty;\mathbb{R}), \\
B(t, \omega, i)\in L_{\mathcal{F}^W}^\infty(0, \infty;\mathbb{R}^m), \\
C(t, \omega, i)\in L_{\mathcal{F}^W}^\infty(0, \infty;\mathbb{R}^n), \\
D(t, \omega, i)\in L_{\mathcal{F}^W}^\infty(0, \infty;\mathbb{R}^{n\times m}), \\
Q(t, \omega, i)\in L_{\mathcal{F}^W}^\infty(0, \infty;\mathbb{R}), \\
R(t, \omega, i)\in L_{\mathcal{F}^W}^\infty(0, \infty;\mathbb{S}^m).
\end{cases}
\end{align*}
\end{assumption}
\begin{assumption}
\label{assump2}
For all $i\in\cM$, $2A(i)+C(i)'C(i)\leq -\rho$, where $\rho> 0$ is a deterministic constant.
\end{assumption}

Under Assumptions \ref{assump1} and \ref{assump2},
there exists a constant $c_1>0$ such that $Q(t,i)\leq c_1$ for all $i\in\cM$. Since $0\in\mathcal{U}$,
\begin{align*}
\inf_{u\in\mathcal{U}}J_\infty(x, i_0, u(\cdot))&\leq J_\infty(x,i_0, 0)\leq c_1x^2\int_0^\infty e^{-\rho t}dt=\frac{c_1}{\rho}x^2.
\end{align*}
So the value function $V(x,i_0)$ of problem \eqref{LQ} is bounded from above. If it is also bounded from below, then the problem is well-posed, namely it has a finite optimal value.

Problem \eqref{LQ} is said to be solvable, if there exists a control $u^*(\cdot)\in\mathcal{U}$ such that
\begin{align*}
-\infty<J_\infty(x, i_0, u^*(\cdot))\leq J_\infty(x, i_0, u(\cdot)), \quad \forall\; u(\cdot)\in\mathcal{U},
\end{align*}
in which case, $u^*(\cdot)$ is called an optimal control for problem \eqref{LQ}, and the optimal value is
\begin{align*}
V(x, i_0):=\inf_{u\in\mathcal{U}}J_\infty(x, i_0, u(\cdot))=J_\infty(x, i_0, u^*(\cdot)).
\end{align*}

\bigskip

\begin{remark}
\label{L2stable}
Under Assumptions \ref{assump1} and \ref{assump2}, for any $u(\cdot)\in L^2_{\mathcal{F}}(0,\infty;\mathbb{R}^m)$, the corresponding state process $X(\cdot)$ of \eqref{state} indeed satisfies both the {$L^2$-stable} condition (see e.g. Definition 2.1 in \cite{SY}):
\begin{align*}
X(\cdot)\in L^2_{\mathcal{F}}(0,\infty;\mathbb{R}),
\end{align*}
and the mean-square stabilizing condition (see e.g. Definition 2.1 in \cite{CZ} or \cite{LZR}):
\[\lim\limits_{T\rightarrow\infty}\E[X(T)^2]=0.\]

In fact, applying It\^{o}'s lemma to $X(\cdot)^2$, we have (the argument $(t,\alpha_t)$ are suppressed for notation simplicity),
\begin{align*}
X(T)^2&=x^2+\int_0^T\Big((2A+C'C)X^2+2X(B+D'C)'u+u'D'Du\Big)dt+\int_0^T2 X(C X+ Du)'dW.
\end{align*}
Let
\[\tau_n=\inf\Big\{T\geq 0: \int_0^T |X(C X+ Du)|^2dt >n\Big\}.\]
Because $X(\cdot)$ is continuous, it is locally bounded almost surely. Hence, $\{\tau_n\}$ is a non-decreasing sequence of stopping times such that $\lim_{n\to\infty}\tau_n=\infty$ and
\begin{align*}
\E[X(T\wedge \tau_n)^2]&=x^2+\E\int_0^{T\wedge \tau_n}\Big((2A+C'C)X^2+2X(B+D'C)'u+u'D'Du\Big)dt.
\end{align*}
By Assumptions \ref{assump1}, \ref{assump2} and the elementary inequality $2ab\leq \frac{\rho}{2}a^2+\frac{2}{\rho}b^2$, we have
\begin{align*}
0\leq \E[X(T\wedge \tau_n)^2]
&\leq x^2+\E\int_0^{T\wedge \tau_n}\Big(-\rho X^2+ \frac{\rho}{2}X^2+\frac{2}{\rho}|(B+D'C)'u|^2+u'D'Du\Big)dt\nonumber\\
&\leq x^2+\E\int_0^{T\wedge \tau_n}\Big(-\frac{\rho}{2}X^2+c|u|^2\Big)dt,
\end{align*}
for some deterministic constant $c>0$, so
\begin{align*}
\frac{\rho}{2}\E\int_0^{T\wedge \tau_n}X^2dt\leq x^2+c\E\int_0^{T\wedge \tau_n}|u|^2dt\leq x^2+c\E\int_0^\infty|u|^2dt.
\end{align*}
First passing $n\rightarrow\infty$ and then passing $T\rightarrow\infty$, by the monotone convergence theorem, we proved that $X(\cdot)$ satisfies the {$L^2$-stable} condition $X(\cdot)\in L^2_{\mathcal{F}}(0,\infty;\mathbb{R})$. As a consequence,
there exists a deterministic sequence $\{T_i\}$ with $\lim\limits_{i\rightarrow\infty}T_i=\infty$ such that
$$\lim_{i\rightarrow\infty}\E[X(T_i)^2]=0.$$
Fix $i$.
Let
\[\theta_n=\inf\Big\{S\geq T_i: \int_{T_i}^S |X(t)(C(t)X(t)+ D(t)u(t))|^2dt >n\Big\}.\]
Then $\{\theta_n\}$ is a non-decreasing sequence of stopping times such that $\lim_{n\rightarrow\infty}\theta_n=\infty$.
For $T>T_i$, applying It\^{o}'s lemma to $X(\cdot)^2$ on $[T_i, T\wedge \theta_n]$, and using a similar argument as above, we get
\begin{align*}
\E[X(T\wedge \theta_n)^2]
&\leq \E[X(T_i)^2]+\E\int_{T_i}^{T\wedge \theta_n}\Big(-\frac{\rho}{2}X^2+c|u|^2\Big)dt\leq \E[X(T_i)^2]+\E\int_{T_i}^{\infty} c|u|^2 dt.
\end{align*}
Passing $n\rightarrow\infty$ and using Fatou's lemma, we obtain
\begin{align*}
\E[X(T)^2] \leq \E[X(T_i)^2]+\E\int_{T_i}^{\infty} c|u|^2 dt,
\end{align*}
which gives
\begin{align*}
\limsup_{T\rightarrow\infty} \E[X(T)^2] \leq \E[X(T_i)^2]+\E\int_{T_i}^{\infty} c|u|^2 dt.
\end{align*}
Because the right hand side convergences to 0 as $i\to \infty$, we conclude
\begin{align*}
\lim_{T\rightarrow\infty} \E[X(T)^2]=0.
\end{align*}
\end{remark}

\section{The extended stochastic Riccati equations}
To tackle problem \eqref{LQ}, we first introduce two related $\ell$-dimensional BSDEs on infinite time horizon.

For $\Lambda\in\mathbb{R}^n$ and $P\in\mathbb{R}$ with $PD(t, i)'D(t, i)+R(t, i)>0$, set
\begin{align*}
H_1(t, \omega, P, \Lambda, i)&=\inf_{v\in\Gamma}\big[v'(PD(t, i)'D(t, i)+R(t, i))v
+2v'(PB(t, i)+PD(t, i)'C(t, i)+D(t, i)'\Lambda)\big], \\
H_2(t, \omega, P, \Lambda, i)&=\inf_{v\in\Gamma}\big[v'(PD(t, i)'D(t, i)+R(t, i))v
-2v'(PB(t, i)+PD(t, i)'C(t, i)+D(t, i)'\Lambda)\big].
\end{align*}
Because $PD(t, i)'D(t, i)+R(t, i)$ is positive definite, $H_{1}$ and $H_{2}$ are well-defined, that is, $\R$-valued. Clearly, they are non-positive as $0\in\Gamma$.


We introduce the following two $\ell$-dimensional ESREs on infinite time horizon (remind that the arguments $t$ and $\omega$ are suppressed):
\begin{align}
\label{P1}
\begin{cases}
dP_1(i)=-\Big[(2A(i)+C(i)'C(i))P_1(i)+2C(i)'\Lambda_1(i)+Q(i)\\
\qquad\qquad\qquad+H_1(P_1(i), \Lambda_1(i), i)+\sum\limits_{j=1}^{\ell}q_{ij}P_1(j)\Big]dt+\Lambda_1(i)'dW, \\
R(i)+P_1(i)D(i)'D(i)>0, \ \mbox{ for all $i\in\cM$};
\end{cases}
\end{align}
and
\begin{align}
\label{P2}
\begin{cases}
dP_2(i)=-\Big[(2A(i)+C(i)'C(i))P_2(i)+2C(i)'\Lambda_2(i)+Q(i)\\
\qquad\qquad\qquad+\ H_2(P_2(i), \Lambda_2(i), i)+\sum\limits_{j=1}^{\ell}q_{ij}P_2(j)\Big]dt+\Lambda_2(i)'dW, \\
R(i)+P_2(i)D(i)'D(i)>0, \ \mbox{ for all $i\in\cM$.}
\end{cases}
\end{align}
The equations \eqref{P1} and \eqref{P2} are referred to as the ESREs. By a solution $(P_{1}(i),\Lambda_{1}(i))_{i\in\cM}$ to \eqref{P1}, we mean $(P_{1}(i),\Lambda_{1}(i))_{i\in\cM}$ satisfying \eqref{P1} in an arbitrary finite time horizon and $(P_1(i),\Lambda_1(i))\in L^\infty_{\mathcal{F}^W}(0,\infty; \mathbb{R})\times L^{2, \;\mathrm{loc}}_{\mathcal {F}^W}(0, \infty;\mathbb{R}^n)$ for all $i\in\cM$. Furthermore, a solution of \eqref{P1} is called nonnegative (respectively uniformly positive) if $P_1(i)\geq 0$ (respectively $P_1(i)\geq c$ for some constant $c>0$) for all $i\in\cM$.
Solutions for ESRE \eqref{P2} are defined similarly.

\begin{remark}
If $\Gamma$ is symmetric, namely, $-v\in\Gamma$ whenever $v\in\Gamma$, then
\( H_1(P, \Lambda, i)=H_2(P, \Lambda, i).\) In particular, if there is no control constraint, i.e. $\Gamma=\mathbb{R}^m$, then they are both equal to
\begin{align*}
-[PB(i)+D(i)'(PC(i)+\Lambda)]'(R(i)+PD(i)'D(i))^{-1}[PB(i)+D(i)'(PC(i)+\Lambda)].
\end{align*}
Furthermore, \eqref{P1} and \eqref{P2} reduce to the same $\ell$-dimensional stochastic Riccati equation on infinite time horizon:
\begin{align*}
\begin{cases}
dP(i)=-\Big[(2A(i)+C(i)'C(i))P(i)+2C(i)'\Lambda(i)+Q(i)\\
\qquad\qquad\qquad-[PB(i)+D(i)'(PC(i)+\Lambda)]'(R(i)+PD(i)'D(i))^{-1}
[PB(i)+D(i)'(PC(i)+\Lambda)]\\
\qquad\qquad\qquad+\sum\limits_{j=1}^{\ell}q_{ij}P(j)\Big]dt+\Lambda(i)'dW, \\
R(i)+P(i)D(i)'D(i)>0, \ \mbox{ for all $i\in\cM$}.
\end{cases}
\end{align*}

\end{remark}

\begin{remark}
If $\ell=1$ (namely, there is no regime switching), then \eqref{P1} and \eqref{P2} become, respectively, the following ESREs on infinite time horizon:
\begin{align*}
\begin{cases}
dP_1=-\Big[(2A+C'C)P_1+2C'\Lambda_1+Q+H_1(P_1, \Lambda_1)\Big]dt+\Lambda_1'dW, \\
R+P_1D'D>0;
\end{cases}
\end{align*}
and
\begin{align*}
\begin{cases}
dP_2=-\Big[(2A+C'C)P_2+2C'\Lambda_2+Q+\ H_2(P_2, \Lambda_2)\Big]dt+\Lambda_2'dW, \\
R+P_2D'D>0.
\end{cases}
\end{align*}
For the case that $m=n=1, \ \Gamma=\mathbb{R}^+$ and $R$ is uniformly positive, these two equations are studied in \cite{PZ}.
\end{remark}

\begin{remark}
If all the coefficients $A, B, C, D, Q, R$ are deterministic continuous and bounded functions of $(t,i)$, then \eqref{P1} and \eqref{P2} have deterministic solutions, where
$\Lambda_{1}(i)=\Lambda_{2}(i)=0$ for all $i\in\cM$, and $P_1(i)$, $P_2(i)$ satisfy, respectively, the following systems of ordinary differential equations (ODE) on infinite time horizon:
\begin{align*}
\begin{cases}
dP_1(i)=-\Big[(2A(i)+C(i)'C(i))P_1(i)+Q(i)+H_1(P_1(i), 0, i)+\sum\limits_{j=1}^{\ell}q_{ij}P_1(j)\Big]dt, \\
R(i)+P_1(i)D(i)'D(i)>0, \ \mbox{ for all $i\in\cM$};
\end{cases}
\end{align*}
and
\begin{align*}
\begin{cases}
dP_2(i)=-\Big[(2A(i)+C(i)'C(i))P_2(i)+Q(i)+\ H_2(P_2(i), 0, i)+\sum\limits_{j=1}^{\ell}q_{ij}P_2(j)\Big]dt, \\
R(i)+P_2(i)D(i)'D(i)>0, \ \mbox{ for all $i\in\cM$.}
\end{cases}
\end{align*}
Furthermore, if all the coefficients $A, B, C, D, Q, R$ only depend on $i$, then
they are also equivalent to the following systems of deterministic algebraic Riccati equations:
\begin{align*}
\begin{cases}
(2A(i)+C(i)'C(i))P_1(i)+Q(i)+H_1(P_1(i), 0, i)+\sum\limits_{j=1}^{\ell}q_{ij}P_1(j)=0, \\
R(i)+P_1(i)D(i)'D(i)>0, \ \mbox{ for all $i\in\cM$};
\end{cases}
\end{align*}
and
\begin{align*}
\begin{cases}
(2A(i)+C(i)'C(i))P_2(i)+Q(i)+ H_2(P_2(i), 0, i)+\sum\limits_{j=1}^{\ell}q_{ij}P_2(j)=0, \\
R(i)+P_2(i)D(i)'D(i)>0, \ \mbox{ for all $i\in\cM$}.
\end{cases}
\end{align*}
\end{remark}

\begin{remark}
If all the coefficients $A, B, C, D, Q, R$ are constant matrices, and $\Gamma=\mathbb{R}^m$, then
$P_{1}(1)=P_2(1)=\cdot\cdot\cdot=P_{1}(\ell)=P_{2}(\ell)$,
$\Lambda_{1}(i)=\Lambda_{2}(i)=0$ for all $i\in\cM$. Both \eqref{P1} and \eqref{P2} become
\begin{align*}
\begin{cases}
dP=-\Big[(2A+C'C)P+Q-P^2(B'+C'D)(R+PD'D)^{-1}(B+D'C)\Big]dt,\\
R+PD'D>0,
\end{cases}
\end{align*}
which is equivalent to
\[
(2A+C'C)P+Q-P^2(B'+C'D)(R+PD'D)^{-1}(B+D'C)=0.
\]
Furthermore, if $2A+C'C<0$, $B\neq0$, $D\equiv0$, $Q>0$, $R>0$, then \eqref{P1} has two solutions
\begin{align*}
\tilde P=\frac{-(2A+C'C)+\sqrt{(2A+C'C)^2+4B'R^{-1}BQ}}{2(2A+C'C)}<0;\\
\hat P=\frac{-(2A+C'C)-\sqrt{(2A+C'C)^2+4B'R^{-1}BQ}}{2(2A+C'C)}>0.
\end{align*}
So the solution of \eqref{P1} as well as that of \eqref{P2} may not be unique. However, we will show the nonnegative solutions to them are unique in the following two sections.
\end{remark}

In the sequel, we will study the solvability of ESREs \eqref{P1} and \eqref{P2} as well as problem \eqref{LQ} in two cases separately: (1) standard case, in which $R(i)$ is uniformly positive definite; (2) singular case, in which $R(i)$ is positive semidefinite. Here ``singular" means that the control weight matrix $R(i)$ in the cost functional \eqref{costfunctional} could be probably a singular matrix.

\section{Standard case}
In this section, we will study the solvability of ESREs \eqref{P1} and \eqref{P2} as well as problem \eqref{LQ} under the following assumption.
\begin{assumption}[Standard case]\label{assump3}
$Q(i)\geq0$ and $R(i)\geq\delta I_{m\times m}$ for all $i\in\cM$, where $\delta$ is some positive constant.
\end{assumption}
Under Assumption \ref{assump3}, clearly we have $J_\infty(x, i_0, u(\cdot))\geq 0$, for all $(x, i_0, u)\in\R\times\cM\times\mathcal{U}$. Therefore, problem \eqref{LQ} is well-posed with a nonnegative optimal value.
\subsection{Solvability of the extended stochastic Riccati equations}
Our first result on the solvability of ESREs \eqref{P1} and \eqref{P2} is given as follows.
\begin{theorem}[Standard case]
\label{existence}
Under Assumptions \ref{assump1}, \ref{assump2} and \ref{assump3}, there exists a nonnegative solution $(P_{1}(i),\Lambda_{1}(i))_{i\in\cM}$ to \eqref{P1} (respectively, $(P_{2}(i),\Lambda_{2}(i))_{i\in\cM})$ to \eqref{P2}).
\end{theorem}
In the rest of this subsection, we prove Theorem \ref{existence}. We only prove the existence of a solution to \eqref{P1}, and the argument for \eqref{P2} is similar. Our method is first to approximate \eqref{P1} by a sequence of BSDEs on finite time horizon, and then pass to the limit.

For any $N>0$, consider the following system of BSDEs on $[0,N]$:
\begin{align}
\label{P1N}
\begin{cases}
dP_{1,N}(i)=-\Big[(2A(i)+C(i)'C(i))P_{1,N}(i)+2C(i)'\Lambda_{1,N}(i)+Q(i)\\
\qquad\qquad\qquad+H_1(P_{1,N}(i), \Lambda_{1,N}(i), i)+\sum\limits_{j=1}^{\ell}q_{ij}P_{1,N}(j)\Big]dt+\Lambda_{1,N}(i)'dW, \\
P_{1,N}(N,i)=0,\\
R(i)+P_{1,N}(i)D(i)'D(i)>0, \ \mbox{ for all $i\in\cM$}.
\end{cases}
\end{align}
By Theorem 3.5 in \cite{HSX}, under Assumptions \ref{assump1}, \ref{assump2} and \ref{assump3}, BSDE \eqref{P1N} admits a unique solution $(P_{1,N}(i),\Lambda_{1,N}(i))_{i\in\cM}$, such that $(P_{1,N}(i), \ \Lambda_{1,N}(i))\in L^\infty_{\mathcal{F}^W}(0, N; \mathbb {R})\times \BMO$ and $P_{1,N}(i)\geq 0$ for all $i\in\cM$. We simply extend $(P_{1,N}(i),\Lambda_{1,N}(i))$ to $[0,+\infty)$ by setting $P_{1,N}(t,i)=0, \ \Lambda_{1,N}(t,i)=0$ for $t>N$. Clearly the extended solution is nonnegative.

In order to pass to the limit, we need to show the solution is upper bounded, uniformly in $N$, and monotone increasing with respect to $N$. They are contained in the subsequent two lemmas.
\begin{lemma}
\label{upperbound}
Under Assumptions \ref{assump1}, \ref{assump2} and \ref{assump3},
we have $P_{1,N}(i)\leq \frac{c_1}{\rho}$ for all $i\in \cM$, where $c_1>0$ is a constant such that $Q(i)\leq c_1$ for all $i\in\cM$.
\end{lemma}
\pf
Notice that $H_1(t,P,\Lambda)\leq0$,
from the proof of Theorem 3.5 in \cite{HSX}, we have
\begin{align*}
P_{1,N}(t, i)\leq \frac{c_1}{\rho}(1-e^{-\rho(N-t)})\leq \frac{c_1}{\rho}.
\end{align*}
So the solution of \eqref{P1N} is upper bounded uniformly in $N$.
\eof

\begin{lemma}
\label{monotone}
Under Assumptions \ref{assump1}, \ref{assump2} and \ref{assump3}, we have $P_{1,N}(i)\leq P_{1,M}(i)$, if $N\leq M$.
\end{lemma}
\pf
Let $a>0$ be sufficiently small such that $R(i)-aD(i)'D(i)> 0$ for all $i\in\cM$. For every $i\in\cM$, define
\[
(U_M(t,i),V_M(t,i))=\left(\ln\left(P_{1,M}(t,i)+a\right),\;\frac{\Lambda_{1,M}(t,i)}{P_{1,M}(t,i)+a}\right),
\]
and
\[
(U_N(t,i),V_N(t,i))=\left(\ln\left(P_{1,N}(t,i)+a\right),\;\frac{\Lambda_{1,N}(t,i)}{P_{1,N}(t,i)+a}\right).
\]
Then $\ln a\leq U_M(t,i)\leq \ln(c_1/\rho+a)$ by Lemma \ref{upperbound}. By It\^{o}'s lemma, $(U_M(t,i),V_M(t,i))_{i\in\cM}$ satisfies the following $\ell$-dimensional BSDE:
\begin{align*}
\begin{cases}
dU_M(i)=-\Big[(2A(i)+C(i)'C(i))(1-ae^{-U_M(i)})+2C(i)'V_M(i)+Q(i)e^{-U_M(i)}\\
\qquad\qquad\qquad+\tilde H(U_M(i), V_M(i), i)+\frac{1}{2}V_M(i)'V_M(i)+\sum\limits_{j=1}^{\ell}q_{ij}e^{U_M(j)-U_M(i)}\Big]dt+V_M(i)'dW, \\
U_M(M, i)=\ln a, \ \mbox{ for all $i\in\cM$, }
\end{cases}
\end{align*}
where
\[\tilde H(U, V, i)=\inf_{v\in\Gamma}\Big[v'((1-ae^{-U})D(i)'D(i)+R(i)e^{-U})v+2v'((1-ae^{-U})(B(i)+D(i)'C(i))+D(i)'V)\Big].\]
Set
\[
(\bar U(t,i), \ \bar V(t,i))=(U_M(t,i)-U_N(t,i), \ V_M(t,i)-V_N(t,i)).
\]
Noticing that $U_M(t,i)=\ln a, V_M(t,i)=0$ for $t>M$, we have
\begin{align*}
\bar U(t,i)&=\int_t^M\Big[(Q(i)-2aA(i)-aC(i)'C(i))(e^{-U_M(i)}-e^{-U_N(i)})+2C(i)'\bar V(i)+\frac{Q(i)}{a}I_{\{s\geq N\}}\\
&\qquad\qquad+\tilde H(U_M(i),V_M(i),i)-\tilde H(U_N(i),V_N(i),i)+(V_M(i)+V_N(i))'\bar V(i)\\
&\qquad\qquad+\sum\limits_{j=1}^{\ell}q_{ij}(e^{U_M(j)-U_M(i)}-e^{U_N(j)-U_N(i)})\Big]ds-\int_t^M\bar V(i)'dW.
\end{align*}
Furthermore, by It\^{o}'s lemma,
\begin{align*}
(\bar U(t,i)^-)^2&=\int_t^M\Big[-2\bar U(i)^-(Q(i)-2aA(i)-aC(i)'C(i))(e^{-U_M(i)}-e^{-U_N(i)})\\
&\qquad\qquad-4\bar U(i)^-C(i)'\bar V(i)-2\bar U(i)^-\frac{Q(i)}{a}I_{\{s\geq N\}}\\
&\qquad\qquad-2\bar U(i)^-(\tilde H(U_M(i),V_M(i),i)-\tilde H(U_N(i),V_N(i),i))\\
&\qquad\qquad-2\bar U(i)^-(V_M(i)+V_N(i))'\bar V(i)-\bar V(i)^2I_{\{\bar U(i)<0\}}\\
&\qquad\qquad-2\bar U(i)^-\sum\limits_{j=1}^{\ell}q_{ij}(e^{U_M(j)-U_M(i)}-e^{U_N(j)-U_N(i)})\Big]ds+\int_t^M2\bar U(i)^-\bar V(i)dW.
\end{align*}
Notice that $-2\bar U(i)^-\frac{Q(i)}{a}I_{\{s\geq N\}}\leq 0$,
using similar method as Theorem 3.5 in \cite{HSX}, we can prove $\bar U(t,i)^-=0$ for a.e. $t\in[0,M]$ and all $i\in\cM$.
This completes the proof.
\eof

We are now ready to prove Theorem \ref{existence} by passing to the limit.

$\textbf{Proof of Theorem \ref{existence}}$.
For any fixed $T>0$, we choose $N>T$ in \eqref{P1N}.
From Lemma \ref{monotone} and Lemma \ref{upperbound}, $P_N(i)$ is non-decreasing w.r.t. $N$ and has a uniform upper bound $\frac{c_1}{\rho}$. Therefore we can define an $\{\mathcal{F}^W_t\}_{t\geq 0}$-predictable process $P_1^T(t,i)$ by \[P_1^T(t,i):=\lim\limits_{N\rightarrow\infty}P_{1,N}(t,i), \mbox{ for all $i\in\cM$.}\]
Note that $P_{1}^T(i)$ does not depend on $T$, and we take this notation just for convenience.

Note that $P_{1,N}(i)\leq\frac{c_1}{\rho}$ and
\begin{align*}
|H_1(P_{1,N}(i), \Lambda_{1,N}(i), i)|&\leq \frac{1}{\delta}|P_{1,N}(i)B(i)+P_{1,N}(i)D(i)'C(i)+D(i)'\Lambda_{1,N}(i)|^2\\
&\leq\frac{c_2}{\delta}(P_{1,N}(i)^2+|\Lambda_{1,N}(i)|^2),
\end{align*}
for some sufficiently large $c_2$, so we can regard $(P_{1,N}(i),\Lambda_{1,N}(i))$ as the solution of a scalar-valued quadratic BSDE on $[0,T]$ for each $i\in\cM$. Thus by Proposition 2.4 in \cite{Ko}, there exists a process $\Lambda_1^T(i)\in L^2_{\mathcal{F}^W}(0,T;\mathbb{R}^n)$ such that $(P_1^T(i),\Lambda_1^T(i))$ is a solution to the $i^{\mbox{th}}$ equation in ESRE \eqref{P1} on the time interval $[0,T]$ with terminal value $P_1^T(T,i)$.

We need to show that $\Lambda_1^T(i)$ does not depend on $T$. For any $T_1> T$, there exists a process $\Lambda_1^{T_1}(i)\in L^2_{\mathcal{F}^W}(0,T_1;\mathbb{R}^n)$ such that $(P_1^{T_1}(i),\Lambda_1^{T_1}(i))$ is a solution to the $i^{\mbox{th}}$ equation in ESRE \eqref{P1} on the time interval $[0,T_1]$ with terminal value $P_1^{T_1}(T_1,i)$. Notice that $P_1^{T_1}(T,i)=\lim\limits_{N\rightarrow\infty}P_{1,N}(T,i)=P_1^{T}(T,i)$, by the uniqueness of the solution of the $i^{\mbox{th}}$ equation in ESRE \eqref{P1} with terminal value $P_1^{T_1}(T,i)=P_1^{T}(T,i)$ (see Theorem 3.5 of \cite{HSX}), we get $\Lambda_1^{T}(t,i)=\Lambda_1^{T_1}(t,i), \ t\in[0,T]$.

This shows that $(P_{1}(i),\Lambda_{1}(i))_{i\in\cM}$ forms a \emph{nonnegative} solution to ESRE \eqref{P1}. This completes the proof.
\eof

The following proposition shows that the solution $(P_{1}(i),\Lambda_{1}(i))_{i\in\cM}$ of \eqref{P1} constructed in Theorem \ref{existence} is actually uniformly positive if $Q(i)$ is uniformly positive for all $i\in\cM$.
\begin{proposition}
If Assumptions \ref{assump1}, \ref{assump2}, \ref{assump3} hold and $Q(i)\geq\delta$ for all $i\in\cM$, where $\delta$ is some positive constant. Then the solution $(P_{1}(i),\Lambda_{1}(i))_{i\in\cM}$ of \eqref{P1} constructed in the proof of Theorem \ref{existence} is uniformly positive.
\end{proposition}
\pf
Let $c_3>0$ be any constant such that
\[
2A(i)+C(i)'C(i)+q_{ii}-\frac{2c_1}{\delta\rho}|B(i)+D(i)'C(i)|^2>-c_3, \ \mbox{for all} \ i\in\cM.
\]
Consider the following $\ell$-dimensional BSDE:
\begin{align}
\label{P1Ndecouple}
\begin{cases}
d\underline P_{1,N}(i)=-\Big[(2A(i)+C(i)'C(i)+q_{ii})\underline P_{1,N}(i)+2C(i)'\underline \Lambda_{1,N}(i)+Q(i)\\
\qquad\qquad\qquad+H_1(\underline P_{1,N}(i),\underline \Lambda_{1,N}(i), i)\Big]dt+\underline \Lambda_{1,N}(i)'dW, \\
\underline P_{1,N}(N,i)=0,\\
R(i)+\underline P_{1,N}(i)D(i)'D(i)>0, \ \mbox{ for all $i\in\cM$}.
\end{cases}
\end{align}
This is a decoupled system of BSDE.
From Theorem 4.1 and Theorem 5.2 of \cite{HZ}, the $i^{\mbox{th}}$ equation in \eqref{P1Ndecouple} admits a unique, hence maximal solution (see page 565 of \cite{Ko} for its definition) $(\underline P_{1,N}(i), \ \underline\Lambda_{1,N}(i))\in L^\infty_{\mathcal{F}^W}(0, N; \mathbb {R})\times \BMO$, and $\underline P_{1,N}(i)\geq0$ for all $i\in\cM$. From the proof of Theorem 3.5 of \cite{HSX}, the solution $(P_{1,N}(i), \ \Lambda_{1,N}(i))_{i\in\cM}$ of \eqref{P1N} could be approximated by the solutions of a sequence of BSDEs with Lipschitz generators. Thus we can use comparison theorem for multidimensional BSDEs (see for example Lemma 3.4 of \cite{HSX}) and then pass to the limit to get $P_{1,N}(i)\geq\underline P_{1,N}(i)$ for all $i\in\cM$.

Let $g:\mathbb{R}^+\rightarrow [0, 1]$ be a smooth truncation function satisfying $g(x)=1$ for $x\in[0, \frac{c_1}{\rho}]$, and $g(x)=0$ for $x\in[2\frac{c_1}{\rho}, +\infty)$. Notice that
$\frac{c_1}{\rho}\geq P_{1,N}(i)\geq\underline P_{1,N}(i)$, so $(\underline P_{1,N}(i), \ \underline\Lambda_{1,N}(i))$ is still a solution of the $i^{\mbox{th}}$ equation in BSDE \eqref{P1Ndecouple} when $H_1(P,\Lambda,i)$ is replaced by $H_1(P,\Lambda,i)g(P)$ in the generator.

Notice that for $P=\underline P_{1,N}(i), \ \Lambda=\underline\Lambda_{1,N}(i))$, we have
\begin{align*}
&\quad\;H_1(P, \Lambda, i)g(P)\\
&=\inf_{v\in\Gamma}\big[v'(PD(i)'D(t, i)+R(i))v
+2v'(PB(i)+PD(i)'C(i)+D(i)'\Lambda)\big]g(P)\\
&\geq\inf_{v\in\Gamma}\big[\delta|v|^2
+2v'(PB(i)+PD(i)'C(i)+D(i)'\Lambda)\big]g(P)\\
&\geq\inf_{v\in\mathbb{R}^m}\big[\delta|v|^2
+2v'(PB(i)+PD(i)'C(i)+D(i)'\Lambda)\big]g(P)\\
&=-\frac{1}{\delta}|PB(i)+PD(i)'C(i)+D(i)'\Lambda|^2g(P)\\
&=-\frac{P^2}{\delta}|B(i)+D(i)'C(i)|^2g(P)
-\frac{P}{\delta}\big(B(i)+D(i)'C(i)\big)'D(i)'\Lambda g(P)-\frac{1}{\delta}|D(i)'\Lambda|^2 g(P)\\
&\geq-\frac{2c_1P}{\delta\rho}|B(i)+D(i)'C(i)|^2
-\frac{P}{\delta}\big(B(i)+D(i)'C(i)\big)'D(i)'\Lambda g(P)-\frac{1}{\delta}|D(i)'\Lambda|^2 g(P).
\end{align*}
The following BSDE
\begin{align*}
\begin{cases}
dP=-\Big[-c_3P+\delta+2C(i)'\Lambda
-\frac{P}{\delta}\big(B(i)+D(i)'C(i)\big)'D(i)'\Lambda g(P)-\frac{1}{\delta}|D(i)'\Lambda|^2 g(P)\Big]dt+\Lambda'dW,\\
P(N)=0,
\end{cases}
\end{align*}
admits a solution $(\frac{\delta}{c_3}(1-e^{-c_3(N-t)}),0)$. Then the maximal solution argument (Theorem 2.3 of \cite{Ko}) gives
\[
\underline P_{1,N}(t,i)\geq\frac{\delta}{c_3}(1-e^{-c_3(N-t)}).
\]
Therefore, for any $t\in[0,\infty)$ and $i\in \cM$, we get
\[
P_1(t,i)=\lim_{N\rightarrow\infty}P_{1,N}(t,i)\geq\lim_{N\rightarrow\infty}\underline P_{1,N}(t,i)
\geq\lim_{N\rightarrow\infty}\frac{\delta}{c_3}(1-e^{-c_3(N-t)})=\frac{\delta}{c_3}.
\]
This completes the proof.
\eof

We now turn to the study of the original stochastic LQ problem \eqref{LQ}.

\subsection{Solution to problem \eqref{LQ}}
For $PD(t, i)'D(t, i)+R(t, i)>0$,
there exists $\tilde R_P(t,i)\in L^\infty_{\mathcal{F}^W}(0,\infty;\mathbb{R}^{m\times m})$ such that
\begin{align*}
PD(t, i)'D(t, i)+R(t, i)=\tilde R_P(t,i)'\tilde R_P(t,i).
\end{align*}
Denote
$\mbox{Proj}_{\tilde R_P\Gamma}(\cdot)$ be the projection mapping from $\mathbb{R}^m$ to the closed cone $\tilde R_P\Gamma$ under the Euclidean norm.\footnote{We remark that the cone $\tilde R_P\Gamma$ depends on $P$.}
The set is not empty and may contain more than one point.
And any point in $\mbox{Proj}_{\tilde R_P\Gamma}\big(-\tilde R_P(PD'D+R)^{-1}(PB+PD'C+D'\Lambda)\big)$ attains the infimum in the definition of $H_1$ because (remind that the arguments $t$, $i$ and $\omega$ are suppressed in $B,\ C, \ D, \ R, \ \tilde R_P$)
\begin{align*}
H_1(t, \omega, P, \Lambda, i)&=\inf_{v\in\Gamma}\big[v'(PD'D+R)v
+2v'(PB+PD'C+D'\Lambda)\big]\\
&=\inf_{v\in\Gamma}\big[v'\tilde R_P'\tilde R_Pv+2v'(PB+PD'C+D'\Lambda)\big]\\
&=\inf_{v\in\Gamma}\Big|\tilde R_P v+\tilde R_P(PD'D+R)^{-1}(PB+PD'C+D'\Lambda)\Big|^2\\
&\quad\qquad-(PB+PD'C+D'\Lambda)'(PD'D+R)^{-1}(PB+PD'C+D'\Lambda).
\end{align*}
Similarly, any point in $\mbox{Proj}_{\tilde R_P\Gamma}\big(\tilde R_P(PD'D+R)^{-1}(PB+PD'C+D'\Lambda)\big)$ attains the infimum in the definition of $H_2$.
From a measurable selection theorem (see e.g. Lemma 11 of \cite{HIM}), we know there exist two $\{\mathcal{F}^W_t\}_{t\geq 0}$-predictable processes
$\hat v_1(t, \omega, P, \Lambda, i)$ and $\hat v_2(t, \omega, P, \Lambda, i)$ that satisfy
\begin{align}\label{hatu}
&\hat v_1(t, \omega, P, \Lambda, i)\in\mbox{Proj}_{\tilde R_P\Gamma}\big(-\tilde R_P(PD'D+R)^{-1}(PB+PD'C+D'\Lambda)\big),\nonumber\\
&\hat v_2(t, \omega, P, \Lambda, i)\in\mbox{Proj}_{\tilde R_P\Gamma}\big(\tilde R_P(PD'D+R)^{-1}(PB+PD'C+D'\Lambda)\big).
\end{align}
We will give an optimal control to problem \eqref{LQ} by $\hat v_1$, $\hat v_2$.

\begin{theorem}
\label{optivalue}
Suppose Assumptions \ref{assump1}, \ref{assump2} and \ref{assump3} hold. Let $(P_{1}(i),\Lambda_{1}(i))_{i\in\cM}$ (resp. $(P_{2}(i),\Lambda_{2}(i))_{i\in\cM}$) be nonnegative solutions to \eqref{P1} (resp. $\eqref{P2}$), and $\hat v_1$, $\hat v_2$ defined in \eqref{hatu}. Then problem \eqref{LQ} admits an optimal control, as a feedback function of the time $t$, the state $X$, and the market regime $i$,
\begin{align}
\label{opticon}
u^*(t, X, i)=\hat v_1(t, P_1(t, i), \Lambda_1(t, i), i)X^++\hat v_2(t, P_2(t, i), \Lambda_2(t, i), i)X^-.
\end{align}
Moreover, the corresponding optimal value is
\begin{align*}
V(x, i_0)=P_1(0, i_0)(x^+)^2+P_2(0, i_0)(x^-)^2.
\end{align*}
\end{theorem}
\pf
For any $T>0$, we consider the following optimal control problems on finite time horizon $[0,T]$:
\begin{align}
\begin{cases}
\mathrm{Minimize} &\ J_{T}(x, i_0, u(\cdot))\\
\mbox{subject to} &\ (X(\cdot), u(\cdot)) \mbox{ admissible for} \ \eqref{state},
\end{cases}
\label{LQT}%
\end{align}
where
\begin{align*}
J_{T}(x, i_0, u(\cdot))&:=\mathbb{E}\bigg\{\int_0^T\Big(Q(t, \alpha_t)X(t)^2+u(t)'R(t, \alpha_t)u(t)\Big)dt\nonumber\\
&\qquad\qquad+P_{1}(T,\alpha_T)(X(T)^+)^2+P_{2}(T,\alpha_T)(X(T)^-)^2\bigg\}.
\end{align*}
According to Theorem 4.2 in \cite{HSX},
\[
u^*(t,X,i)=\hat v_1(t, P_{1}(t, i), \Lambda_{1}(t, i), i)X^++\hat v_2(t, P_{2}(t, i), \Lambda_{2}(t, i), i)X^-
\]
is an optimal feedback control of problem \eqref{LQT} with the optimal value
\begin{align*}
\inf_{u\in\mathcal{U}} J_{T}(x,i_0,u(\cdot))= J_{T}(x,i_0,u^*(\cdot))=P_{1}(0,i_0)(x^+)^2+P_{2}(0,i_0)(x^-)^2.
\end{align*}
Therefore, for any admissible pare $(X(\cdot),u(\cdot))$ for \eqref{state}, we have
\begin{align*}
P_{1}(0,i_0)(x^+)^2+P_{2}(0,i_0)(x^-)^2&\leq J_{T}(x,i_0,u(\cdot))\\
&=\mathbb{E}\bigg\{\int_0^T\Big(Q(t, \alpha_t)X(t)^2+u(t)'R(t, \alpha_t)u(t)\Big)dt\\
&\qquad\qquad+P_{1}(T,\alpha_T)(X(T)^+)^2+P_{2}(T,\alpha_T)(X(T)^-)^2\bigg\}\\
&\leq\mathbb{E}\bigg\{\int_0^\infty\Big(Q(t, \alpha_t)X(t)^2+u(t)'R(t, \alpha_t)u(t)\Big)dt\\
&\qquad\qquad+P_{1}(T,\alpha_T)(X(T)^+)^2+P_{2}(T,\alpha_T)(X(T)^-)^2\bigg\}.
\end{align*}
Thanks to Remark \ref{L2stable},
$$\lim_{T\rightarrow\infty}\E[X(T)^2]=0.$$
Since $P_1$ and $P_2$ are uniformly bounded, we have
\begin{align*}
P_{1}(0,i_0)(x^+)^2+P_{2}(0,i_0)(x^-)^2&\leq\lim_{T\rightarrow\infty}\mathbb{E}\bigg\{\int_0^\infty\Big(Q(t, \alpha_t)X(t)^2+u(t)'R(t, \alpha_t)u(t)\Big)dt\\
&\qquad\qquad+P_{1}(T,\alpha_{T})(X(T)^+)^2+P_{2}(T,\alpha_{T})(X(T)^-)^2\bigg\}\\
&=\mathbb{E}\int_0^\infty\Big(Q(t, \alpha_t)X(t)^2+u(t)'R(t, \alpha_t)u(t)\Big)dt.
\end{align*}

Taking infimum over $u\in\mathcal{U}$, we obtain
\begin{align}\label{leq}
P_{1}(0,i_0)(x^+)^2+P_{2}(0,i_0)(x^-)^2\leq \inf_{u\in\mathcal{U}}J_\infty (x,i_0,u(\cdot)).
\end{align}

On the other hand, let $X^*(t)$ be the associated solution of \eqref{state} with $u$ replaced by $u^*$. Then
\begin{align*}
P_{1}(0,i_0)(x^+)^2+P_{2}(0,i_0)(x^-)^2&=\inf_{u\in\mathcal{U}} J_{T}(x,i_0,u(\cdot))= J_{T}(x,i_0,u^*(\cdot))\\
&=\mathbb{E}\bigg\{\int_0^T\Big(Q(t, \alpha_t)X^*(t)^2+u^*(t)'R(t, \alpha_t)u^*(t)\Big)dt\\
&\qquad\qquad+P_{1}(T,\alpha_T)(X^*(T)^+)^2+P_{2}(T,\alpha_T)(X^*(T)^-)^2\bigg\}\\
&\geq\mathbb{E}\int_0^T\Big(Q(t, \alpha_t)X^*(t)^2+u^*(t)'R(t, \alpha_t)u^*(t)\Big)dt,
\end{align*}
where the last inequality is due to the nonnegativity of $P_1$ and $P_2$.
Note that $X^{*}$ and $u^*$ are independent of $T$.
Letting $T\rightarrow\infty$ in above, by the monotone convergence theorem, we get
\begin{align*}
P_1(0,i_0)(x^+)^2+P_2(0,i_0)(x^-)^2&\geq \mathbb{E}\int_0^\infty\Big(Q(t, \alpha_t) X^*(t)^2+ u^*(t)'R(t, \alpha_t) u^*(t)\Big)dt.
\end{align*}
By Assumption \ref{assump3}, $Q\geq 0$ and $R\geq\delta I_m$, so the right hand side is
\begin{align*}
\geq \delta\mathbb{E}\int_0^\infty|u^*(t)|^2dt.
\end{align*}
Hence $u^*\in L^2_{\mathcal{F}}(0,\infty;\mathbb{R}^m)$ and $u^*\in\mathcal{U}$. This further implies
\begin{align}
\label{geq}
P_1(0,i_0)(x^+)^2+P_2(0,i_0)(x^-)^2&\geq \mathbb{E}\int_0^\infty\Big(Q(t, \alpha_t) X^*(t)^2+ u^*(t)'R(t, \alpha_t) u^*(t)\Big)dt\nonumber\\
&\geq\inf_{u\in\mathcal{U}}J_{\infty}(x,i_0,u(\cdot)).
\end{align}

Combining \eqref{leq} and \eqref{geq}, we conclude that
\[
\inf_{u\in\mathcal{U}}J_\infty (x,i_0,u(\cdot))=P_{1}(0,i_0)(x^+)^2+P_{2}(0,i_0)(x^-)^2= J_{\infty}(x,i_0,u^*(\cdot)),
\]
and consequently, $u^*(t,X,\alpha_t)$ is an optimal feedback control for \eqref{LQ}.
\eof

As a byproduct of Theorem \ref{optivalue}, we have the following uniqueness of the solutions for ESREs \eqref{P1} and \eqref{P2}.
\begin{theorem}\label{uniqestan}
Under Assumptions \ref{assump1}, \ref{assump2} and \ref{assump3}, each of the ESREs \eqref{P1} and \eqref{P2} admits at most one nonnegative solution.
\end{theorem}
\pf
Consider an stochastic LQ control problem on $[s,\infty)$, with $s\geq0$, where the system dynamics is \eqref{state} with initial time $s$, initial state $x_s\in L^2_{\mathcal{F}_s}(\Omega;\mathbb{R})$ and initial regime $\alpha_s=i$, and the cost functional is
\begin{align*}
J_{s,\infty}(x_s,i,u):=\mathbb{E}\Big\{\int_s^\infty\Big(Q(t, \alpha_t)X(t)^2+u(t)'R(t, \alpha_t)u(t)\Big)dt\;\Big|\;X_s=x_s, \alpha_s=i\Big\}.
\end{align*}
Let $(P_{1}(i),\Lambda_{1}(i))_{i\in\cM}$ and $(P_{2}(i),\Lambda_{2}(i))_{i\in\cM}$ be nonnegative solutions to \eqref{P1} and $\eqref{P2}$ respectively. Then going through the same analysis as in the proof of Theorem \ref{optivalue}, we deduce that the optimal cost is
\begin{align*}
V(x_s,i):=\inf_{u \ \mathrm{admissible}}J_{s,\infty}(x_s,i,u)=P_{1}(s,i)(x_s^+)^2+P_{2}(s,i)(x_s^-)^2,
\end{align*}
which clearly implies the uniqueness.
\eof
\section{Singular case}
In this section, we will study the solvability of ESREs \eqref{P1} and \eqref{P2} as well as problem \eqref{LQ} in a singular case.
\begin{assumption}[Singular case]
\label{assump4}
$R(i)\geq0$, and there exists a constant $\delta>0$ such that $D(i)'D(i)\geq \delta I_{m\times m}$ and $Q(i)\geq\delta$ for all $i\in\cM$.
\end{assumption}
Under Assumption \ref{assump4}, clearly $J_\infty(x, i_0, u(\cdot))\geq 0$, for all $(x, i_0, u)\in\R\times\cM\times\mathcal{U}$, so problem \eqref{LQ} is well-posed.
\subsection{Solvability of the extended stochastic Riccati equations}
Same as before, we first study the solvability of ESREs \eqref{P1} and \eqref{P2}.
\begin{theorem}[Singular case]
\label{singular}
Under Assumptions \ref{assump1}, \ref{assump2} and \ref{assump4}, there exists a uniformly positive solution $(P_{1}(i),\Lambda_{1}(i))_{i\in\cM}$ to \eqref{P1} (respectively, $(P_{2}(i),\Lambda_{2}(i))_{i\in\cM})$ to \eqref{P2}).
\end{theorem}

In the rest of this subsection, we prove Theorem \ref{singular}. We only prove the existence of a solution to \eqref{P1}, and the argument for \eqref{P2} is similar.

For any $a>0$, consider the following system of BSDEs on $[0,N]$,
\begin{align}
\label{P1Na}
\begin{cases}
dP_{1,N}^{a}(i)=-\Big[(2A(i)+C(i)'C(i))P_{1,N}^a(i)+2C(i)'\Lambda_{1,N}^a(i)+Q(i)\\
\qquad\qquad\qquad+H_1^a(P_{1,N}^a(i), \Lambda_{1,N}^a(i), i)+\sum\limits_{j=1}^{\ell}q_{ij}P_{1,N}^a(j)\Big]dt+\Lambda_{1,N}^a(i)'dW, \\
P_{1,N}^a(N,i)=0,\\
aI_m+R(i)+P_{1,N}^a(i)D(i)'D(i)>0, \ \mbox{ for all $i\in\cM$},
\end{cases}
\end{align}
where
\begin{align*}
H_1^a(t, \omega, P, \Lambda, i)&=\inf_{v\in\Gamma}\big[v'(PD(t, i)'D(t, i)+R(t, i)+aI_m)v
+2v'(PB(t, i)+PD(t, i)'C(t, i)+D(t, i)'\Lambda)\big].
\end{align*}
By Theorem 3.5 in \cite{HSX}, under Assumptions \ref{assump1}, \ref{assump2} and \ref{assump4}, BSDE \eqref{P1Na} admits a unique solution $(P_{1,N}^a(i),\Lambda^a_{1,N}(i))_{i\in\cM}$, such that $(P^a_{1,N}(i), \ \Lambda^a_{1,N}(i))\in L^\infty_{\mathcal{F}^W}(0, N; \mathbb {R})\times \BMO$ and $P^a_{1,N}(i)\geq 0, $ for all $i\in\cM$. Same as before, we extend $(P^a_{1,N}(i),\Lambda^a_{1,N}(i))$ to $[0,+\infty)$ by setting $P^a_{1,N}(t,i)=0, \ \Lambda^a_{1,N}(t,i)=0$ for $t>N$.

The following lemma studies the monotonicity of the solution of BSDE \eqref{P1Na} with respect to $a$.
The proof is similar to Lemma \ref{opticon}, so we left the details to interested readers.
\begin{lemma}
\label{monotonea}
Under Assumptions \ref{assump1}, \ref{assump2} and \ref{assump4}, we have
$P_{1,N}^{a_1}(t,i)\geq P_{1,N}^{a_2}(t,i)\geq 0$, if $a_1\geq a_2>0$.
\end{lemma}

Our next result gives a lower bound, uniformly in $a$, for the the solution of BSDEs \eqref{P1Na}.
\begin{lemma}
\label{lbsin}
Suppose Assumptions \ref{assump1}, \ref{assump2} and \ref{assump4} hold.
Then
\begin{align*}
P^a_{1,N}( i)\geq \frac{\delta}{c_4}(1-e^{-c_4(N-t)}), \ \mbox{for a.e. $t\in[0,N]$ and all $i\in\cM$,}
\end{align*}
where $c_4>0$ is any constant satisfying
\[
2A(i)+C(i)'C(i)+q_{ii}-\frac{1}{\delta}|B(i)+D(i)'C(i)|^2\geq -c_4, \ \mbox{for all} \ i\in\cM.
\]
\end{lemma}
\pf
Consider the following BSDE
\begin{align}
\label{underlineP}
\begin{cases}
d\underline P(i)=-\Big[(2A(i)+C(i)'C(i)+q_{ii})\underline P(i)+2C(i)'\underline\Lambda(i)+Q(i)\\
\qquad\qquad\qquad+H_1^a(\underline P(i), \underline\Lambda(i), i)\Big]dt+\underline\Lambda(i)'dW, \\
\underline P(i,N)=0,\\
aI_m+R(i)+\underline P(i)D(i)'D(i)>0, \ \mbox{ for all $i\in\cM$}.
\end{cases}
\end{align}
This is a decoupled system of BSDE.
By Theorem 4.1 and Theorem 5.2 of \cite{HZ}, there exists a unique solution $(\underline P(i),\underline\Lambda(i))\in L^\infty_{\mathcal{F}^W}(0, N; \mathbb {R})\times \BMO$ with $\underline P(i)\geq 0 $ for all $i\in\cM$.

Under Assumptions \ref{assump1}, \ref{assump2} and \ref{assump4}, for $P\geq 0$,
we have
\begin{align*}
H_1^a(P, \Lambda, i)&=\inf_{v\in\Gamma}\big[v'(PD(i)'D(i)+R( i)+aI_m)v+2v'(PB(i)+PD(i)'C(i)+D(i)'\Lambda)\big]\\
&\geq \inf_{v\in\Gamma}\big[(\delta P+a) |v|^2+2v'(PB(i)+PD(i)'C(i)+D(i)'\Lambda)\big]\\
&\geq \inf_{v\in\mathbb{R}^m}\big[(\delta P+a)|v|^2+2v'(PB(i)+PD(i)'C(i)+D(i)'\Lambda)\big]\\
&=-\frac{P^2}{\delta P+a}|B(i)+D(i)'C(i)|^2
-\frac{2P}{\delta P+a}\big(B(i)+D(i)'C(i)\big)'D(i)'\Lambda-\frac{1}{\delta P+a}|D(i)'\Lambda|^2\\
&\geq-\frac{P}{\delta }|B(i)+D(i)'C(i)|^2
-\frac{2P}{\delta P+a}\big(B(i)+D(i)'C(i)\big)'D(i)'\Lambda-\frac{1}{\delta P+a}|D(i)'\Lambda|^2.
\end{align*}
Consider the following BSDE:
\begin{align*}
\begin{cases}
dP=-\Big[-c_4P+\delta+2C(i)'\Lambda
-\frac{2P}{\delta P+a}\big(B(i)+D(i)'C(i)\big)'D(i)'\Lambda\\
\qquad\qquad-\frac{1}{\delta P+a}|D(i)'\Lambda|^2\Big]dt+\Lambda'dW,\\
P(N)=0.
\end{cases}
\end{align*}
Obviously it admits a solution$(\frac{\delta}{c_4}(1-e^{-c_4(N-t)}),0)$. Then from the maximal solution argument (Theorem 2.3 of \cite{Ko}), we know that the solution $(\underline P(i), \underline\Lambda(i))$ of \eqref{underlineP} satisfies
\begin{align*}
\underline P(t, i)\geq \frac{\delta}{c_4}(1-e^{-c_4(N-t)}).
\end{align*}
Then, similar to Theorem \ref{existence}, we have
\begin{align*}
P^a_{1,N}(t, i)\geq\underline P(t, i)\geq \frac{\delta}{c_4}(1-e^{-c_4(N-t)}), \ \mbox{for a.e. $t\in[0,N]$ and all $i\in\cM$}.
\end{align*}
The proof is finished.
\eof

We are ready to prove Theorem \ref{singular}.

\textbf{Proof of Theorem \ref{singular}}.
By Theorem \ref{existence} and \ref{uniqestan}, the following BSDEs on infinite horizon admit a unique nonnegative solution $(P_1^a(i),\Lambda^a_1(i))_{i\in\cM}$:
\begin{align}
\label{P1a}
\begin{cases}
dP_{1}^{a}(i)=-\Big[(2A(i)+C(i)'C(i))P_{1}^a(i)+2C(i)'\Lambda_{1}^a(i)+Q(i)\\
\qquad\qquad\qquad+H_1^a(P_{1}^a(i), \Lambda_{1}^a(i), i)+\sum\limits_{j=1}^{\ell}q_{ij}P_{1}^a(j)\Big]dt+\Lambda_{1}^a(i)'dW, \\
aI_{m}+R(i)+P_{1}^a(i)D(i)'D(i)>0, \ \mbox{ for all $i\in\cM$}.
\end{cases}
\end{align}
And $P^a_1(t,i)=\lim\limits_{N\rightarrow\infty}P_{1,N}^a(t,i)$.
From Lemma \ref{lbsin}, we immediately get, for any $t>0$ $$P^a_1(t,i)=\lim\limits_{N\rightarrow\infty}P_{1,N}^a(t,i)
\geq\lim\limits_{N\rightarrow\infty}\frac{\delta}{c_4}(1-e^{-c_4(N-t)})=\frac{\delta}{c_4}.$$
For any $T>0$, $(P_1^a(i),\Lambda^a_1(i))_{i\in\cM}$ is the unique nonnegative solution of the following system of BSDEs on $[0,T]$:
\begin{align}
\begin{cases}
dP_{1}^{a}(i)=-\Big[(2A(i)+C(i)'C(i))P_{1}^{a}(i)+2C(i)'\Lambda_{1}^{a}(i)+Q(i)\\
\qquad\qquad\qquad+H_1^a(P_{1}^{a}(i), \Lambda_{1}^{a}(i), i)+\sum\limits_{j=1}^{\ell}q_{ij}P_{1}^{a}(j)\Big]dt+\Lambda_{1}^{a}(i)'dW, \\
P_1^{a}(T,i)=P_1^a(T,i),\\
aI_{m}+R(i)+P_{1}^{a}(i)D(i)'D(i)>0, \ \mbox{ for all $i\in\cM$}.
\end{cases}
\end{align}
Recalling that $P_1^a(i)\geq\frac{\delta}{c_4}$ and by Lemma \ref{monotonea},
we can set $P_1(t,i)=\lim\limits_{a\rightarrow0}P_1^{a}(t,i)$ for all $i\in\cM$. Note that we can regard
$(P_1^a(i),\Lambda_1^a(i))$ as the solution of a scalar-valued quadratic BSDE on $[0,T]$ for each $i\in\cM$. Let $a\rightarrow0$, by Proposition 2.4 in \cite{Ko}, there exists a process $\Lambda_1(i)\in L^2_{\mathcal{F}^W}(0,T;\mathbb{R}^n)$ such that $(P_1(i),\Lambda_1(i))$ is a solution to ESRE \eqref{P1} on the time interval $[0,T]$ with terminal value $P_1(T,i)$. Similar argument as in Theorem \ref{existence} yields that $\Lambda_1(i)$ is independent of $T$. Thus $(P_1(i),\Lambda_1(i))_{i\in\cM}$ is a uniformly positive solution to \eqref{P1}.
\eof

%

We now ready to solve the stochastic LQ problem \eqref{LQ} in the singular case.
\subsection{Solutions to problem \eqref{LQ}}

\begin{theorem}
\label{optivaluesin}
Suppose Assumptions \ref{assump1}, \ref{assump2} and \ref{assump4} hold. Let $(P_{1}(i),\Lambda_{1}(i))_{i\in\cM}$ and $(P_{2}(i),\Lambda_{2}(i))_{i\in\cM}$ be uniformly positive solutions to \eqref{P1} and $\eqref{P2}$, respectively. Then problem \eqref{LQ} admits an optimal control, as a feedback function of the time $t$, the state $X$, and the market regime $i$,
\begin{align}
\label{opticonsin}
u^*(t, X, i)=\hat v_1(t, P_1(t, i), \Lambda_1(t, i), i)X^++\hat v_2(t, P_2(t, i), \Lambda_2(t, i), i)X^-.
\end{align}
Moreover, the corresponding optimal value is
\begin{align}
\label{VeqP}
V(x,i_0)=P_1(0, i_0)(x^+)^2+P_2(0, i_0)(x^-)^2.
\end{align}
\end{theorem}
\pf
We use the same notation as Theorem \ref{optivalue}. By similar analysis as in the proof of Theorem \ref{optivalue}, we can deduce that
\begin{align}\label{leqsin}
P_{1}(0,i_0)(x^+)^2+P_{2}(0,i_0)(x^-)^2\leq \inf_{u\in\mathcal{U}}J_\infty (x,i_0,u(\cdot)).
\end{align}

New we show $u^*$ defined in \eqref{opticonsin} is an optimal feedback control of problem \eqref{LQ}.
Let $X^*(t)$ be the associated solution of \eqref{state} with $u$ replaced by $u^*$. By Theorem 4.2 in \cite{HSX},
\begin{align*}
P_{1}(0,i_0)(x^+)^2+P_{2}(0,i_0)(x^-)^2&=\inf_{u\in\mathcal{U}} J_{T}(x,i_0,u(\cdot))= J_{T}(x,i_0,u^*(\cdot))\\
&=\mathbb{E}\bigg\{\int_0^T\Big(Q(t, \alpha_t)X^*(t)^2+u^*(t)'R(t, \alpha_t)u^*(t)\Big)dt\\
&\qquad\qquad+P_{1}(T,\alpha_T)(X^*(T)^+)^2+P_{2}(T,\alpha_T)(X^*(T)^-)^2\bigg\}\\
&\geq\mathbb{E}\int_0^T\Big(Q(t, \alpha_t)X^*(t)^2+u^*(t)'R(t, \alpha_t)u^*(t)\Big)dt.
\end{align*}
By monotone convergence theorem, we have
\begin{align}
\label{upper}
\mathbb{E}\int_0^\infty\Big(Q(t, \alpha_t)X^*(t)^2+u^*(t)'R(t, \alpha_t)u^*(t)\Big)dt\leq P_{1}(0,i_0)(x^+)^2+P_{2}(0,i_0)(x^-)^2.
\end{align}
This together with \eqref{leqsin} clearly implies $u^*$ is an optimal control of problem \eqref{LQ}, if we can show $u^*\in L^2_{\mathcal{F}}(0,\infty;\mathbb{R}^m)$. To show $u^*\in L^2_{\mathcal{F}}(0,\infty;\mathbb{R}^m)$, we first notice by Assumption \ref{assump4} and \eqref{upper},
\begin{align*}
\delta\mathbb{E}\int_0^\infty X^*(t)^2dt\leq
\mathbb{E}\int_0^\infty Q(t, \alpha_t)X^*(t)^2dt\leq P_{1}(0,i_0)(x^+)^2+P_{2}(0,i_0)(x^-)^2,
\end{align*}
so $X^*\in L^2_{\mathcal{F}}(0,\infty;\mathbb{R})$.
Applying It\^{o}'s lemma to $X^*(t)^2$,
\begin{align}\label{hatxeq}
X^*(t)^2-x^2&=\int_0^t\Big((2A+C'C)(X^*)^2+2 X^*(B'+C'D) u^*+|D u^*|^2\Big)ds\nonumber\\
&\qquad\qquad+\int_0^t2 X^*(C X^*+ D(u^*))'dW.
\end{align}
For $n=1,2,...$, set
\[
\tau_n=\inf\left\{t>0:\int_0^t|2 X^*(C' X^*+ (u^*)'D')|^2ds\geq n\right\}\wedge n.
\]
Then $\tau_n\uparrow \infty$ as $n\rightarrow\infty$ and it follows from \eqref{hatxeq} that
\begin{align*}
\E \big( X^*(t\wedge\tau_n)^2\big)&=x^2+\E\int_0^{t\wedge\tau_n}
\Big((2A+C'C) (X^*)^2+2X^*(B'+C'D) u^*+|Du^*|^2\Big)ds\\
&\geq x^2+\E\int_0^{t\wedge\tau_n}
\Big((2A+C'C)(X^*)^2-\frac{\delta}{2}| u^*|^2
-\frac{2}{\delta}|B'+C'D|^2(X^*)^2+\delta|u^*|^2\Big)ds\\
&\geq x^2+\E\int_0^{t\wedge\tau_n}
\Big(\frac{\delta}{2}|u^*|^2-\beta (X^*)^2\Big)ds,
\end{align*}
where $\beta>0$ is any constant such that $2A(i)+C(i)'C(i)-\frac{2}{\delta}|B(i)'+C(i)'D(i)|^2\geq-\beta$ for all $i\in\cM$. So
\begin{align*}
\beta\E\int_0^{t\wedge\tau_n} (X^*)^2ds+\E \big(X^*(t\wedge\tau_n)^2\big) \geq x^2+\frac{\delta}{2}\E\int_0^{t\wedge\tau_n}|u^*|^2ds.
\end{align*}
By Lemma 4.3 of \cite{HSX}, we know $u^*\in L^2_{\mathcal{F}}(0,t;\mathbb{R}^m)$ for any $t>0$. So by standard theory of SDE, we have
that $\E \big(\sup_{s\leq t} X^*(s)^{2}\big)<\infty$ for any $t>0$.
Because $X^*(t\wedge\tau_n)^{2} \leq \sup_{s\leq t} X^*(s)^{2}$, letting $n\rightarrow\infty$ in the above inequality and using the dominated convergence theorem and the monotone convergence theorem, we obtain
\begin{align*}
\beta\E\int_0^{t} (X^*)^2ds+\E \big(X^*(t)^2\big) \geq x^2+\frac{\delta}{2}\E\int_0^{t}|u^*|^2ds.
\end{align*}
Let $t\rightarrow\infty$, it follows from the monotone convergence theorem that
\begin{align}
\label{ufinite}
\beta\E\int_0^{\infty} (X^*)^2ds +\lim_{t\rightarrow\infty }\E \big(X^*(t)^2\big)\geq x^2+\frac{\delta}{2}\E\int_0^{\infty}| u^*|^2ds.
\end{align}
Because $X^*\in L^2_{\mathcal{F}}(0,\infty;\mathbb{R})$, the first term on the left hand side is finite. And the second term is 0 from Remark \ref{L2stable}. Hence $u^*\in L^2_{\mathcal{F}}(0,\infty;\mathbb{R}^m)$.

Combining \eqref{leqsin} and \eqref{upper}, we conclude \eqref{VeqP},
and this completes the proof.
\eof

Similar as in the standard case, we have the following uniqueness of the solutions for ESREs \eqref{P1} and \eqref{P2}.
\begin{theorem}\label{uniqesing}
Under Assumptions \ref{assump1}, \ref{assump2} and \ref{assump4}, each of the ESRE \eqref{P1} and \eqref{P2} admits at most one nonnegative solution.
\end{theorem}

We have finished the study of problem \eqref{LQ} in the singular case. In the next section, we apply these results to study a lifetime portfolio selection problem of tracking a given wealth level with regime switching and portfolio constraint.

\section{Tracking a given wealth level}
Consider a financial market consisting of a risk-free asset (the money market
instrument or bond) whose price is $S_{0}$ and $m$ risky securities (the
stocks) whose prices are $S_{1}, \ldots, S_{m}$. Assume $m\leq n$, i.e., the number of risky securities is no more than the dimension of the Brownian motion.
These asset prices are driven by SDEs:
\begin{align*}
\begin{cases}
dS_0(t)=r(t,\alpha_t)S_{0}(t)dt, \\
S_0(0)=s_0,
\end{cases}
\end{align*}
and
\begin{align*}
\begin{cases}
dS_k(t)=S_k(t)\Big(\mu_k(t, \alpha_t)dt+\sum\limits_{j=1}^n\sigma_{kj}(t, \alpha_t)dW_j(t)\Big), \\
S_k(0)=s_k,
\end{cases}
\end{align*}
where $r(t,i)$ is the interest rate process and $\mu_k(t, i)$ and $\sigma_k(t, i):=(\sigma_{k1}(t, i), \ldots, \sigma_{kn}(t, i))$ are the appreciation rate process and volatility rate process of the $k$th risky security corresponding to a market regime $\alpha_t=i$, for every $k=1, \ldots, m$ and $i\in\cM$.

Define the appreciate vector
\begin{align*}
\mu(t, i)=(\mu_1(t, i), \ldots, \mu_m(t, i))',
\end{align*}
and volatility matrix
\begin{align*}
\sigma(t, i)=
\left(
\begin{array}{c}
\sigma_1(t, i)\\
\vdots\\
\sigma_m(t, i)\\
\end{array}
\right)
\equiv (\sigma_{kj}(t, i))_{m\times n}, \ \text{for}\ \text{each} \ i\in\cM.
\end{align*}
In the rest part of this section, we shall assume
$r(\cdot,\cdot,i), \ \mu_k(\cdot, \cdot, i)$, $\sigma_{kj}(\cdot, \cdot, i)\in L^\infty_{\mathcal F^W}(0, \infty;\mathbb R)$, for all $k=1, \ldots, m$, $j=1, \ldots, n$, and $i\in\cM$.

A small investor, whose actions cannot affect the asset prices, will decide at every time
$t\in[0, \infty)$ what amount $\pi_j(t)$ of his wealth to invest in the $j$th risky asset, $j=1, \ldots, m$. The vector process $\pi(\cdot):=(\pi_1(\cdot), \ldots, \pi_m(\cdot))'$ is called a portfolio of the investor. Then the investor's self-financing wealth process $X(\cdot)$ corresponding to a portfolio $\pi(\cdot)$ is a strong solution of the SDE:
\begin{align}
\label{wealth}
\begin{cases}
dX(t)=[r(t,\alpha_t)X(t)+\pi(t)'b(t, \alpha_t)]dt+\pi(t)'\sigma(t, \alpha_t)dW(t), \\
X(0)=x, \ \alpha_0=i_0,
\end{cases}
\end{align}
where $b(t, \alpha_t):=\mu(t, \alpha_t)-r(t,\alpha_t)\mathbf{1}_{m}$ and $\mathbf{1}_{m}$ is the $m$-dimensional vector with all entries being one.


The admissible portfolio set is defined as
\begin{align*}
\mathcal U=\Big\{\pi\in L^2_{\mathcal F}(0, \infty;\mathbb R^m)\;\Big|\; \pi(\cdot)\in\Gamma \Big\},
\end{align*}
where
$\Gamma$ is a given closed cone in $\mathbb R^m$.
For any $\pi\in \mathcal{U}$, the SDE \eqref{wealth} has a unique strong solution.
Economically speaking, no-shorting is allowed in the market if $\Gamma=\R_{+}^{m}$.

\par
For a given wealth level $d\in\mathbb{R}$, the investor's problem is to
\begin{align}
\mathrm{Minimize}&\quad \E\int_0^\infty\bigg(e^{-2\int_0^t\rho(s,\alpha_s)ds}\Big(X(t)-de^{\int_0^tr(s,\alpha_s)ds}\Big)^{2}
+\lambda e^{-2\int_0^t\rho(s,\alpha_s)ds}|\pi_t|^2\bigg)dt%
, \nonumber\\
\mathrm{ s.t.} &\quad
\pi\in \mathcal{U},
\label{optm}%
\end{align}
where $\lambda$ is a real constant, $\rho(\cdot, \cdot, i)\in L^\infty_{\mathcal F^W}(0, \infty;\mathbb R)$ is the discount factor process for every $i\in\cM$. We assume there exists a constant $c_5>0$ such that $\rho(\cdot, \cdot, i)-r(\cdot, \cdot, i)\geq c_5$ for all $i\in\cM$. This ensures that problem \eqref{optm} is well-defined.
\begin{remark}
Economically speaking, problem \eqref{optm} is meaningful only when $x< d$ because nobody would pursuit a target that is not higher than his initial wealth. Mathematically speaking, our argument works for $x\geq d$ as well.
\end{remark}
Besides the above assumptions, we also put the following assumption
in the rest of this section.
\begin{assumption}
Either $\lambda>0$; or
$\lambda=0$ and there exists a constant $\delta>0$ such that $\sigma(i)\sigma(i)'\geq\delta I_{m\times m}$ for all $i\in\cM$.
\end{assumption}

\begin{remark}
It is also possible to consider a model with random coefficient $\lambda$.
\end{remark}
To tackle problem \eqref{optm}, we write $Y(t)=e^{-\int_0^t\rho(s,\alpha_s)ds}\Big(X(t)-de^{\int_0^tr(s,\alpha_s)ds}\Big)$ and $\tilde\pi(t)=e^{-\int_0^t\rho(s,\alpha_s)ds}\pi(t)$. Then by It\^{o}'s lemma,
\begin{align*}
\begin{cases}
dY(t)=[(r(t,\alpha_t)-\rho(t,\alpha_t))Y(t)+\tilde\pi(t)'b(t, \alpha_t)]dt+\tilde\pi(t)'\sigma(t, \alpha_t)dW(t), \\
Y(0)=x-d, \ \alpha_0=i_0,
\end{cases}
\end{align*}
and problem \eqref{optm} becomes a stochastic LQ problem
\begin{align}
\mathrm{Minimize}&\quad \E\int_0^\infty\Big(Y(t)^{2}
+\lambda|\tilde\pi_t|^2\Big)dt%
, \nonumber\\
\mathrm{ s.t.} &\quad
\tilde\pi\in \mathcal{U}.
\end{align}

We now employ the results in Sections 4 and 5 to solve the above problem.
In this case,
\begin{align*}
H_1(t, \omega, P, \Lambda, i)&=\inf_{v\in\Gamma}\big[v'(P\sigma(t, i)\sigma(t, i)'+\lambda I_{m\times m})v
+2v'(Pb(t, i)+\sigma(t, i)\Lambda)\big], \\
H_2(t, \omega, P, \Lambda, i)&=\inf_{v\in\Gamma}\big[v'(P\sigma(t, i)\sigma(t, i)'+\lambda I_{m\times m})v
-2v'(Pb(t, i)+\sigma(t, i)\Lambda)\big].
\end{align*}
And \eqref{P1} and \eqref{P2} become
\begin{align}
\label{P1W}
\begin{cases}
dP_1(i)=-\Big[2(r(i)-\rho(i))P_1(i)+1+H_1(P_1(i), \Lambda_1(i), i)+\sum\limits_{j=1}^{\ell}q_{ij}P_1(j)\Big]dt+\Lambda_1(i)'dW, \\
\lambda I_{m\times m}+P_1(i)\sigma(i)\sigma(i)'>0, \ \mbox{ for all $i\in\cM$};
\end{cases}
\end{align}
and
\begin{align}
\label{P2W}
\begin{cases}
dP_2(i)=-\Big[2(r(i)-\rho(i))P_2(i)+1+\ H_2(P_2(i), \Lambda_2(i), i)+\sum\limits_{j=1}^{\ell}q_{ij}P_2(j)\Big]dt+\Lambda_2(i)'dW, \\
\lambda I_{m\times m}+P_2(i)\sigma(i)\sigma(i)'>0, \ \mbox{ for all $i\in\cM$.}
\end{cases}
\end{align}

From Theorems \ref{optivalue}, \ref{uniqestan}, \ref{optivaluesin} and \ref{uniqesing}, we have the following result for problem \eqref{optm}.
\begin{theorem}
Let $(P_{2}(i),\Lambda_{2}(i))_{i\in\cM}$ be the unique nonnegative solution to $\eqref{P2W}$, and $\hat v_2$ be defined in \eqref{hatu}, i.e.
$\hat v_2(t, \omega, P, \Lambda, i)\in\mbox{Proj}_{\tilde R_P\Gamma}\big(\tilde R_P(P\sigma\sigma'+\lambda I_{m\times m})^{-1}(Pb+\sigma\Lambda)\big).$
If $x\leq d$, then problem \eqref{optm} admits an optimal control, as a feedback function of the time $t$, the state $X$, and the market regime $i$,
\begin{align*}
\pi^*(t, X, i)=\hat v_2(t, P_2(t, i), \Lambda_2(t, i), i)\Big(X(t)-de^{\int_0^tr(s,i)ds}\Big)e^{-\int_0^t\rho(s,i)ds}.
\end{align*}
Moreover, the corresponding optimal value is
\begin{align*}
\inf_{\pi\in\mathcal{U}} \E\int_0^\infty\bigg[e^{-2\int_0^t\rho(s,\alpha_s)ds}\Big(X(t)-de^{\int_0^tr(s,\alpha_s)ds}\Big)^{2}
+\lambda e^{-2\int_0^t\rho(s,\alpha_s)ds}|\pi_t|^2\bigg]dt=P_2(0, i_0)(x-d)^2.
\end{align*}
\end{theorem}

\begin{remark}
If $m=n=1, \ \Gamma=\mathbb{R}, \ \cM=\{1,2\}, \ r=\lambda=0$ and all coefficients are deterministic functions of $i$, then \eqref{P1W} and \eqref{P2W} reduce to
\begin{align*}
\begin{cases}
-2\rho(1)P(1)+1-\frac{b(1)^2}{\sigma(1)^2}P(1)+q_{11}P(1)+q_{12}P(2)=0;\\
-2\rho(2)P(2)+1-\frac{b(2)^2}{\sigma(2)^2}P(2)+q_{21}P(1)+q_{22}P(2)=0;\\
P(i)>0, \ i=1,2.
\end{cases}
\end{align*}
After rearrangement, we get
\begin{align}
\label{system}
\begin{cases}
\left(-2\rho(1)-\frac{b(1)^2}{\sigma(1)^2}+q_{11}\right)P(1)+q_{12}P(2)+1=0;\\
q_{21}P(1)+\left(-2\rho(2)-\frac{b(2)^2}{\sigma(2)^2}+q_{22}\right)P(2)+1=0;\\
P(i)>0, \ i=1,2.
\end{cases}
\end{align}
Noting that $q_{i1}+q_{i2}=0$ for $i=1,2$ and $q_{ij}\geq0$ for $i\neq j$, we know the matrix
\begin{align*}
M=
\begin{pmatrix}
2\rho(1)+\frac{b(1)^2}{\sigma(1)^2}-q_{11} & -q_{12} \\
-q_{21} & 2\rho(2)+\frac{b(2)^2}{\sigma(2)^2}-q_{22}
\end{pmatrix}
\end{align*}
is invertible. Actually,
\begin{align*}
\mathrm{det}(M)&=\left(2\rho(1)+\frac{b(1)^2}{\sigma(1)^2}-q_{11}\right)
\left(2\rho(2)+\frac{b(2)^2}{\sigma(2)^2}-q_{22}\right)-q_{12}q_{21}\\
&=\left(2\rho(1)+\frac{b(1)^2}{\sigma(1)^2}\right)
\left(2\rho(2)+\frac{b(2)^2}{\sigma(2)^2}\right)
+q_{12}\left(2\rho(2)+\frac{b(2)^2}{\sigma(2)^2}\right)
+q_{21}\left(2\rho(1)+\frac{b(1)^2}{\sigma(1)^2}\right)>0.
\end{align*}
Therefore \eqref{system} admits a unique positive solution
\begin{align*}
\begin{pmatrix}
P(1)\\
P(2)
\end{pmatrix}
=\frac{1}{\mathrm{det}(M)}
\begin{pmatrix}
2\rho(2)+\frac{b(2)^2}{\sigma(2)^2}+q_{12}+q_{21} \\
2\rho(1)+\frac{b(1)^2}{\sigma(1)^2}+q_{12}+q_{21}
\end{pmatrix}.
\end{align*}
\end{remark}

\begin{remark}
\label{exam}
If $m=n=1$, $\lambda=0$, $\Gamma=\R_{+}$, and all the coefficients in problem \eqref{optm} are constants. Then \eqref{P2W} admits a unique solution
\begin{align*}
P_2=\frac{1}{2(\rho-r)+\frac{(b^+)^2}{\sigma^2}},
\end{align*}
and the optimal value of problem \eqref{optm} is
\begin{align*}
\inf_{\pi\in \mathcal{U}}\E\int_0^\infty e^{-2\rho t}(X(t)-de^{rt})^2dt=P_2(x-d)^2.
\end{align*}
Moreover, the optimal portfolio of problem \eqref{optm} is
\begin{align*}
\pi^*=-\frac{b^+}{\sigma^2}(X(t)-d)e^{-\rho t}.
\end{align*}
\end{remark}

\begin{remark}
If $\lambda\equiv0$, then
\begin{align*}
\hat v_2(P,\Lambda)&\in\mathrm{Proj}_{\sigma'\Gamma}\left(\sigma'(\sigma\sigma')^{-1}\left(b+\frac{\sigma\Lambda}{P}\right)\right).
\end{align*}
\end{remark}
\begin{remark}
If $\lambda\equiv0$, and $m=n=1$, then
\begin{align*}
\hat v_2(P,\Lambda)=\mathrm{Proj}_{\sigma\Gamma}\left(\frac{b}{\sigma}+\frac{\Lambda}{P}\right),
\end{align*}
is unique.
\end{remark}

\begin{remark}
If $\lambda\equiv0$, and $\rho, \ r, \ b, \ \sigma$ are deterministic functions of $(t,i)$, then $P_2(i)_{i\in\cM}$ is the unique nonnegative solution of ODE
\begin{align*}
\begin{cases}
dP_2(i)=-\Big[2(r(i)-\rho(i))P_2(i)+1+\ H_2(P_2(i), 0, i)+\sum\limits_{j=1}^{\ell}q_{ij}P_2(j)\Big]dt,\\
P_2(i)>0, \ \mbox{ for all $i\in\cM$}.
\end{cases}
\end{align*}
In this case, $\hat v_2$ is independent of $P_2$ and
\begin{align*}
\hat v_2&\in\mathrm{Proj}_{\sigma'\Gamma}\Big(\sigma'(\sigma\sigma')^{-1}b\Big).
\end{align*}
\end{remark}

\bigskip

\section{Concluding remarks}
This paper investigates a stochastic LQ optimal control problem on
infinite time horizon, with regime switching, random coefficients, and cone control constraint. The problem has been completely solved in two different cases, that is, the control weight matrix $R$ in the cost functional being uniformly positive definite or being positive semidefinite. We obtained the optimal
state feedback control and optimal value function by solving two systems of highly nonlinear BSDEs on infinite time horizon. We showed the existence of solutions for the BSDEs by pure BSDE method. But the uniqueness was shown by a verification argument (thus relying on the control problem). Because the solvability of these BSDEs is interesting in its own right from the point view of BSDE theory, we believe it is of great theoretical importance to prove the uniqueness by pure BSDE method in future study. To demonstrate the importance of the theoretical results, we applied the results to solve a lifetime portfolio selection problem of tracking a given wealth level with regime switching and portfolio constraint.

We may consider to extend the present results to several directions, for instance, (1) The constrained LQ control problem when the dimension of state is bigger than 1; (2) The stochastic LQ problem when the control variable is constrained in a convex (but not cone) set such as a bounded interval; (3) Stochastic LQ differential game with control constraints. We hope to address these problems in our future research.

\bigskip
\textbf{Acknowledgment.} The authors would like to thank the anonymous referees for the constructive comments and suggestions, which greatly improve the previous version of the manuscript.

\end{document}